\documentclass[9pt]{amsart}
\usepackage{amsmath}
\usepackage{amsfonts}
\usepackage{graphicx}
\usepackage{amscd}
\usepackage{amsmath,amsfonts,amssymb,amsthm,mathrsfs,xypic}
\newtheorem{thm}{Theorem}[section]

\newtheorem{cor}[thm]{Corollary}
\newtheorem{lem}[thm]{Lemma}
\newtheorem{exm}[thm]{Example}

\newtheorem{con}[thm]{Construction}

\newtheorem{rem}[thm]{Remark}
\newtheorem{defn}[thm]{Definition}
\numberwithin{equation}{section}

\begin{document}

\title[Gorenstein modules]
{On homotopy categories of Gorenstein modules: compact generation and dimensions}
\author[Nan Gao]{Nan  Gao}
\thanks{ e-mail: nangao$\symbol{64}$shu.edu.cn}
\thanks{Supported by the National Natural Science Foundation of China
(Grant No. 11101259).}

\maketitle
\begin{center}
Department of Mathematics, \ \ Shanghai  University\\
Shanghai 200444, P. R. China
\end{center}

\vskip 5pt

\begin{abstract} \ Let $A$ be a virtually Gorenstein  algebra of finite CM-type. We establish a duality between the  subcategory of compact objects in  the homotopy category  of Gorenstein projective left $A$-modules and the bounded Gorenstein derived category of finitely generated right $A$-modules. Let $R$ be a two-sided noetherian ring such that the subcategory of Gorenstein flat modules $R\mbox{-}\mathcal{GF}$ is closed under direct products. We show that the inclusion  $K(R\mbox{-}\mathcal{GF})\to K(R\mbox{-}{\rm Mod})$ of homotopy categories admits a right adjoint. We introduce the notion of Gorenstein representation dimension for an algebra of finite CM-type, and establish relations among the dimension of its relative Auslander algebra, Gorenstein representation dimension,
the dimension of the bounded Gorenstein derived category, and the dimension of the bounded homotopy category of its Gorenstein projective modules.
\vskip 10pt


\vskip 20pt

\end{abstract}

\section{Introduction}

\vskip 10pt

Gorenstein projective modules and algebras of finite CM-type are of interest in the representation theory of
algebras, Gorenstein homological algebra, and in the theory of singularity categories
(see e.g. [AM], [EJ],  [CFH], [Buc], [Be]).

\vskip 10pt

Triangulated categories, especially, derived categories, introduced by Grothendieck and Verdier([Ver]), have been playing an increasingly important role in various areas of mathematics, including representation theory, algebraic geometry, and mathematical physics.
In the last decade, some of the progress has been on Brown representability in homotopy theory (see e.g. [CN], [BS2]).
There are several theorems telling us about the existence and uniqueness of model structures for large classes of triangulated categories ([LO], [S]).
To be able to use Brown representability Theorem, a triangulated category needs a basic property of compact generation ([N2]).
Later, Krause ([K1]) and J${\o}$rgensen ([J]) have established relations between the compact generation of a triangulated category and  the existence of  dualizing complexes, respectively.

\vskip 10pt

A major topic of current interest is the compact generation of the homotopy category of projective modules, $K(R\mbox{-}\mathcal{P})$, of a ring $R$.
J${\o}$rgensen ([J]) has shown for  any reasonably nice ring,  $K(R\mbox{-}\mathcal{P})$ is compactly generated, also he has established  a duality  between its subcategory of compact objects and the bounded derived category of finitely presented right modules. This was generalized by Neeman in [N4]  to arbitrary rings. He has proved that $K(R\mbox{-}\mathcal{P})$ is always $\aleph_{1}$-compactly generated for any ring $R$.
\vskip 10pt

These results raise  questions on the homotopy category of Gorenstein projective modules:

\vskip 10pt

$(1)$ \ When is the homotopy category of Gorenstein projective modules  compactly generated? What is the subcategory of compact objects?
How is it  related to the Gorenstein derived category of the corresponding ring introduced by Gao and Zhang ([GZ])?

\vskip 10pt

$(2)$ \ In this case, what is the relation between the subcategory of compact objects in the homotopy category of Gorenstein projective modules and the subcategory of compact objects in the homotopy category of projective modules?

\vskip 10pt

We completely answer the question $(1)$ and $(2)$ for virtually Gorenstein algebras of finite CM-type (Theorem 3.2 and 3.7).

\vskip 10pt

It is known that projective modules and flat modules are same over artin algebras. However, it is not true in general. For an arbitrary ring $R$, Neeman
([N4]) has shown the inclusion of the homotopy category of flat $R$-modules, $K(R\mbox{-}\mathcal{F})$, into  $K(R\mbox{-}{\rm Mod})$ has a right adjoint. In [K2] a criterion has been given for the  existence of  right approximations in cocomplete additive categories which is used to construct adjoint functors in homotopy categories.

 \vskip 10pt

These results raise questions on the homotopy category of  Gorenstein flat modules:

\vskip 10pt

$(3)$ \ Can we establish a pair of adjoint functors between the homotopy category of whose Gorenstein flat modules and  the homotopy category of
a ring?

\vskip 10pt

We completely answer the question $(3)$ for two-sided noetherian rings such that the subcategories of Gorenstein flat modules
are closed under direct products (Theorem 3.9).

\vskip 10pt

The concept of dimension of a triangulated category has been introduced by Rouquier ([Ro]). He  defined and studied
the dimension for a triangulated category in order to prove the representation
dimension of an algebra can be arbitrarily large. He also has shown for an algebra the relations among
the global dimension,  Auslander's representation dimension and the dimension of the bounded derived category.
See [KK] and [O] for more information on this topics.

\vskip 10pt

These results raise questions on the bounded Gorenstein derived category and the homotopy category of Gorenstein projective modules:

\vskip 10pt

$(4)$ \  What is the analogue of Auslander's representation dimension in Gorenstein homological algebra?
Can we establish relations among it, the dimension of the bounded Gorenstein derived category,
the dimension of the bounded homotopy category of its Gorenstein projective modules, and also the dimension of its relative Auslander algebra
for an algebra of finite CM-type?

\vskip 10pt

We provide such analogues and relate them by a chain of inequalities (Definition 4.2, Theorem 4.5 and Theorem 4.6).

\vskip 10pt

Let us end this introduction by mentioning that model structures for Gorenstein derived categories and homotopy categories of Gorenstein projective modules
are investigated, which will appear in the coming paper.

\vskip 20pt

\section{\bf Preliminaries}

\vskip 10pt

In this section we fix notation and recall the main concepts to be used.

\vskip 10pt

Let $A$ be an artin algebra. Denote by $A\mbox{-}{\rm Mod}({\rm resp.}\ A\mbox{-}{\rm mod})$ the category
of (resp. finitely generated) left
$A\mbox{-}$modules, and $A\mbox{-}{\mathcal P}({\rm resp.} \ A\mbox{-}{\rm proj})$ the full subcategory of (resp. finitely generated
)projective $A\mbox{-}$modules. An $A$-module $M$ is said to be
Gorenstein projective in $A\mbox{-}{\rm Mod}({\rm resp.}A\mbox{-}{\rm mod})$, if there is an
exact sequence $P^{\bullet}=\cdots \longrightarrow
P^{-1}\longrightarrow P^{0} \stackrel{d^0}{\longrightarrow}
P^{1}\longrightarrow P^{2}\longrightarrow \cdots$ in $A\mbox{-}{\mathcal P}({\rm resp.} \ A\mbox{-}{\rm proj})$ with
$\rm{Hom}_A(P^{\bullet}, Q)$ exact for any $A\mbox{-}$module $Q$ in
$A\mbox{-}{\mathcal P}({\rm resp.} \ A\mbox{-}{\rm proj})$, such that $M\cong \operatorname{ker}d^0$ (see [EJ]). Denote by
$A\mbox{-}{\mathcal G}\mathcal{P}({\rm resp.} \ A\mbox{-}\mathcal{G}{\rm proj})$ the full subcategory of
Gorenstein projective modules in $A\mbox{-}{\rm Mod}({\rm resp.}\ A\mbox{-}{\rm mod})$. For the notion of Gorenstein injective module we refer to [EJ]. We denote the subcategory of Gorenstein injective $A$-modules by $A\mbox{-}\mathcal{GI}$.

\vskip 10pt

Now recall the notion of Gorenstein flat module. Denote by $A\mbox{-}{\mathcal F}$ the full subcategory of flat $A\mbox{-}$modules. An $A$-module $M$ is said to be Gorenstein flat in $A\mbox{-}{\rm Mod}$, if there is an
exact sequence $F^{\bullet}=\cdots \longrightarrow
F^{-1}\longrightarrow F^{0} \stackrel{d^0}{\longrightarrow}
F^{1}\longrightarrow F^{2}\longrightarrow \cdots$ in $R\mbox{-}{\mathcal F}$ with
$I\otimes_{R}F^{\bullet}$ exact for any injective right $A\mbox{-}$module $I$, such that $M\cong \operatorname{ker}d^0$ (see [EJ]). Denote by
$A\mbox{-}{\mathcal G}\mathcal{F}$ the full subcategory of
Gorenstein flat modules in $A\mbox{-}{\rm Mod}$.

\vskip 10pt

A proper $A\mbox{-}{\mathcal G}proj$-resolution of $A$-module $M$ in $A\mbox{-}{\rm mod}$ is an
exact sequence $E^{\bullet}=\cdots\longrightarrow G_1\longrightarrow
G_0\longrightarrow M\longrightarrow 0$ such that all $G_i\in
A\mbox{-}\mathcal Gproj$, and that
$\operatorname{Hom}_A(G, E^\bullet)$ stays exact for each
$G\in A\mbox{-}\mathcal Gproj$. The second requirement
guarantees the uniqueness of such a resolution in the homotopy
category (the Comparison Theorem; see [EJ], p.169). The coproper $A\mbox{-}{\mathcal G}proj$-resolution is defined dually.

\vskip 10pt

The Gorenstein projective dimension $\mathcal{G} p{\rm
dim}M$ of $M$ in $A\mbox{-}{\rm mod}$ is defined to be the smallest integer $n\geq 0$ such
that there is an exact sequence $0\longrightarrow
G_n\longrightarrow\cdots \longrightarrow G_1\longrightarrow
G_0\longrightarrow M\longrightarrow 0$ with all $G_i\in
A\mbox{-}\mathcal Gproj$, if it exists; and $\mathcal{G} p{\rm dim}M = \infty$ if there is no such exact sequence
of finite length.

\vskip 10pt

A complex $C^\bullet$ of (finitely generated) $A\mbox{-}$modules is $A\mbox{-}{\mathcal G}\mathcal{P}({\rm resp.} \ A\mbox{-}\mathcal{G}proj)$-exact,
if $\operatorname{Hom}_{A}(G, C^\bullet)$ is exact for
any $G\in A\mbox{-}\mathcal{GP}({\rm resp.} \ A\mbox{-}\mathcal{G}proj)$. It is also called  proper
exact for example in [AM]. A chain map \ $f^\bullet: X^\bullet
\longrightarrow Y^\bullet$ is
 an $A\mbox{-}\mathcal Gproj$-quasi-isomorphism, if \
 $\operatorname{Hom}_{ A}(G, f^\bullet)$ is a
 quasi-isomorphism for any $G\in A\mbox{-}\mathcal Gproj$, i.e., there
 are isomorphisms of abelian groups
 ${\rm H}^n\operatorname{Hom}_{ A}(G, f^\bullet): \
 {\rm H}^n\operatorname{Hom}_{ A}(G, X^\bullet)\cong
 {\rm H}^n\operatorname{Hom}_{A}(G, Y^\bullet), \ \forall \
 n\in\Bbb Z, \ \forall \ G\in A\mbox{-}\mathcal Gproj$.

 \vskip 5pt

Following [GZ2], the  (bounded) Gorenstein derived category
$D_{gp}(A\mbox{-}{\rm Mod})({\rm resp.} \  D^b_{gp}(A))$ of $A$ is defined as the Verdier quotient of the (bounded) homotopy category
$K(A\mbox{-}{\rm Mod})$ $({\rm resp.} \  K^b(A\mbox{-}{\rm mod}))$ with respect to the triangulated
subcategory $K_{gpac}(A\mbox{-}{\rm Mod})({\rm resp.} \  $ $K^b_{gpac}(A\mbox{-}{\rm mod}))$ of
$A\mbox{-}{\mathcal GP}({\rm resp.} \ A\mbox{-}\mathcal{G}proj)$-acyclic complexes.

\vskip 10pt

Recall from [BH,\ Be] an artin algebra $A$ is of finite CM-type if there are only finitely many isomorphism classes of finitely generated
indecomposable Gorenstein projective $A$-modules. Suppose $A$ is an artin algebra of finite CM-type, and $G_1, \cdots, G_n$ are
all the pairwise non-isomorphic indecomposable finitely generated
Gorenstein projective $A$-modules, and $G = \bigoplus\limits_{1\le
i\le n}G_i$. Set $\mathcal{G}p(A):={\rm End}_{A}(G)^{\rm op}$, which
we call the relative Auslander algebra of $A$. It is clear that $G$
is an $A\mbox{-}\mathcal{G}p(A)$-bimodule and $\mathcal{G}p(A)$ is
an artin algebra ([ARS, p.27]). Denote by $\mathcal{G}p(A)\mbox{-}{\rm{mod}}$ the category of finitely
generated left $\mathcal{G}p(A)$-modules. Recall from [BR] that an artin algebra $A$ is called virtually Gorenstein if
$A\mbox{-}\mathcal{GP}^{\perp}=^{\perp}A\mbox{-}\mathcal{GI}$.

\vskip 10pt

Let $\mathcal{T}$ be a triangulated category with arbitrary small coproducts. Recall from [N3] that an objects $T\in \mathcal{T}$ is compact
if the functor ${\rm Hom}_{\mathcal{T}}(T, -)$ preserves coproducts. The full subcategory of all compact objects in $\mathcal{T}$ will be denoted $\mathcal{T}^{c}$.

\vskip 20pt

\section{The subcategory $K(A\mbox{-}\mathcal{GP})^{c}$ of compact objects}

\vskip 10pt

In this section we show that if $A$ is a virtually Gorenstein artin algebra of finite CM-type, then the subcategory of compact objects, $K(A\mbox{-}\mathcal{GP})^{c}$, of the homotopy category $K(A\mbox{-}\mathcal{GP})$ is triangular equivalent to the opposite category of the bounded Gorenstein derived category $D^{b}_{gp}(A^{\rm op})$ of $A^{op}$. We also prove that if $R$ is a two-sided noetherian ring such that the category of Gorenstein flat $R$-modules, $R\mbox{-}\mathcal{GF}$,  is closed under direct products, then the inclusion  of the homotopy category of Gorenstein flat $R$-modules, $K(R\mbox{-}\mathcal{GF})$, into the homotopy  category $K(R\mbox{-}{\rm Mod})$ admits a right adjoint.

\vskip 10pt

Let $A$ be a virtually Gorenstein artin algebra of finite CM-type. Then $A^{\rm op}$ is also a virtually Gorenstein artin algebra of finite CM-type.
By [Be, Proposition 4.18] we have that $A\mbox{-}\mathcal{GP}={\rm Add}(A\mbox{-}\mathcal{G}proj)$ and $A^{\rm op}\mbox{-}\mathcal{G}proj$ is contravariantly
finite in $A^{\rm op}\mbox{-}{\rm mod}$. Denote by $( \ \ )^{*}$ the functor ${\rm Hom}_{A}(-,A)$ which dualizes with respect to $A$.

\vskip 10pt

\begin{con} \ Let $M$ be a  finitely generated $A$-module. Then there is a proper $A^{\rm op}\mbox{-}\mathcal{G}proj\mbox{-}$resolution $G^{\bullet}$ with $G^{\bullet}$ in $K^{-,gpb}(A^{\rm op}\mbox{-}\mathcal{G}proj)$ for $M^{*}$ (See [EJ], also [GZ2]). This deduces from [HJ2, Theorem 1.6]  that $G^{\bullet*}$ is a coproper $A\mbox{-}\mathcal{G}proj$-resolution of $M$. Set
$$\mathcal{G}_{2}:=\{G^{\bullet*}[n]|M\in A\mbox{-}{\rm mod}, n\in \mathbb{Z}\}$$
\end{con}

\vskip 10pt

\begin{thm} \  $K(A\mbox{-}\mathcal{GP})$ is a compactly generated triangulated category with $\mathcal{G}_{2}$ as a set of compact generators. Moreover, there is a triangle-equivalence
$$K(A\mbox{-}\mathcal{GP})^{c}\cong D^{b}_{gp}(A^{\rm op})^{\rm op}$$
\end{thm}

\vskip 5pt

\noindent{\bf Proof.} \ By [HJ1, Theorem 3.1] we get that $\mathcal{G}_{2}$ is a set of compact objects in $K(A\mbox{-}\mathcal{GP})$. Also by the proof of [G2, Theorem 2.2] we get that $\mathcal{G}_{2}$ generates $K(A\mbox{-}\mathcal{GP})$. This implies that $K(A\mbox{-}\mathcal{GP})^{c}$ is a full subcategory of $K(A\mbox{-}\mathcal{GP})$ consisting of objects which are finitely built from objects $G^{\bullet}$ in $\mathcal{G}_{2}$, using shifts, distinguished triangles, and direct summands.

\vskip 5pt

Now we set $\mathcal{G}_{2}^{*}:=\{ G^{\bullet*}|G^{\bullet}\in \mathcal{G}_{2}\}$. Denote by $\mathcal{D}$ the full subcategory of $K(A^{\rm op}\mbox{-}\mathcal{GP})$ consisting of objects which are finitely built from objects  in $\mathcal{G}_{2}^{*}$.
Since the canonical chain maps $G^{\bullet}\to G^{\bullet**}$ and $G^{\bullet*}\to G^{\bullet***}$ are isomorphisms, it follows that
$$K(A\mbox{-}\mathcal{GP})^{c}\rightleftarrows \mathcal{D}^{\rm op}$$
are quasi-inverse equivalences of triangulated categories.

\vskip 5pt

Now we claim that $\mathcal{D}$ consists of the objects finitely built from proper $A^{\rm op}\mbox{-}\mathcal{G}proj$-resolutions of all finitely generated $A^{\rm op}$-modules. Then by [GZ2, Theorem 3.6] we get that $\mathcal{D}$ is triangluar equivalent to $D^{b}_{gp}(A^{\rm op})$. This completes the proof.

\vskip 5pt

Suppose that $N$ is a finitely generated $A^{\rm op}$-module, and let
$$E^{\bullet}=\cdots\to E^{-2}\to E^{-1}\to E^{0}\to 0$$
be a proper $A^{\rm op}\mbox{-}\mathcal{G}proj$-resolution of $N$. Now
$$\widetilde{E^{\bullet}}=\cdots\to E^{-3}\to E^{-2}\to 0\to\cdots$$
is a proper $A^{\rm op}\mbox{-}\mathcal{G}proj$-resolution of $Z^{-1}({\widetilde{E^{\bullet}}})$, the $(-1)$st cycle
module of $E^{\bullet}$. Complete $E^{0*}\to E^{-1*}$ with its cokernel,
$$E^{0*}\to E^{-1*}\to M\to 0.$$
Then  $Z^{-1}({\widetilde{E^{\bullet}}})=M^{*}$. This implies $\widetilde{E^{\bullet}}$ is in $\mathcal{G}_{2}^{*}$, also
$E^{\bullet}$ is in $\mathcal{G}_{2}^{*}$. \hfill$\blacksquare$
\vskip 10pt

\begin{cor} \ Let $A$ be a virtually Gorenstein  algebra of finite CM-type. Then

\vskip 5pt

$(1)$ \ there exists the following recollement
$$
 K(A\mbox{-}\mathcal{GP})\begin{smallmatrix}
  \underleftarrow{ \ \ \  l \ \ \  } \\
  \underrightarrow{ \ \ \ i \ \ \ } \\
  \overleftarrow{ \ \ \ r \ \ \ }
\end{smallmatrix}
K(A\mbox{-}{\rm Mod})
\begin{smallmatrix}
  \underleftarrow{ \ \ \  \ \ \ } \\
  \underrightarrow{\ \ \  \ \ \ } \\
  \overleftarrow{ \ \ \  \ \ \ }
\end{smallmatrix}
K(A\mbox{-}{\rm Mod})/ K(A\mbox{-}\mathcal{GP})
$$
In particular, if $A$ is Gorenstein, then we have the recollement of the form
$$
 K(A\mbox{-}\mathcal{GP})\begin{smallmatrix}
  \underleftarrow{ \ \ \  l \ \ \  } \\
  \underrightarrow{ \ \ \ i \ \ \ } \\
  \overleftarrow{ \ \ \ r \ \ \ }
\end{smallmatrix}
K(A\mbox{-}{\rm Mod})
\begin{smallmatrix}
  \underleftarrow{ \ \ \  \ \ \ } \\
  \underrightarrow{\ \ \  \ \ \ } \\
  \overleftarrow{ \ \ \  \ \ \ }
\end{smallmatrix}
K_{gpac}(A\mbox{-}{\rm Mod})
$$

$(2)$ \ there exists the following right recollement
 $$
 K_{gpac}(A\mbox{-}{\rm Mod})\begin{smallmatrix}
 \underrightarrow{ \ \ \ i \ \ \ } \\
  \overleftarrow{ \ \ \ r \ \ \ }
\end{smallmatrix}
K(A\mbox{-}{\rm Mod})
\begin{smallmatrix}
  \underrightarrow{\ \ \  \ \ \ } \\
  \overleftarrow{ \ \ \  \ \ \ }
\end{smallmatrix}
D_{gp}(A\mbox{-}{\rm Mod})
$$
In particular, $D_{gp}(A\mbox{-}{\rm Mod})$ has small Hom-sets.
In this case, if $A$ is also Gorenstein, then we have the recollement of the form
$$
 K_{gpac}(A\mbox{-}{\rm Mod})\begin{smallmatrix}
  \underleftarrow{ \ \ \  l \ \ \  } \\
  \underrightarrow{ \ \ \ i \ \ \ } \\
  \overleftarrow{ \ \ \ r \ \ \ }
\end{smallmatrix}
K(A\mbox{-}{\rm Mod})
\begin{smallmatrix}
  \underleftarrow{ \ \ \  \ \ \ } \\
  \underrightarrow{\ \ \  \ \ \ } \\
  \overleftarrow{ \ \ \  \ \ \ }
\end{smallmatrix}
K(A\mbox{-}\mathcal{GP})
$$
\end{cor}

\vskip 5pt

\noindent{\bf Proof.} \ $(1)$ \ By Theorem 3.2 we get that $K(A\mbox{-}\mathcal{GP})$ is  compactly generated. Since the inclusion $i$ naturally  preserves  coproducts and products, it follows from  [N2, Theorem 4.1] that $i$  admits  a right adjoint $r$, also a left adjoint $l$.
So by [Mi, Theorem 2.2] we have the following recollement
$$
 K(A\mbox{-}\mathcal{GP})\begin{smallmatrix}
  \underleftarrow{ \ \ \  l \ \ \  } \\
  \underrightarrow{ \ \ \ i \ \ \ } \\
  \overleftarrow{ \ \ \ r \ \ \ }
\end{smallmatrix}
K(A\mbox{-}{\rm Mod})
\begin{smallmatrix}
  \underleftarrow{ \ \ \  \ \ \ } \\
  \underrightarrow{\ \ \  \ \ \ } \\
  \overleftarrow{ \ \ \  \ \ \ }
\end{smallmatrix}
K(A\mbox{-}{\rm Mod})/ K(A\mbox{-}\mathcal{GP})
$$
If $A$ is Gorenstein, then by [G1, Theorem 2.7] we have a triangle-equivalence
$$K(A\mbox{-}{\rm Mod})/ K(A\mbox{-}\mathcal{GP})\cong K_{gpac}(A\mbox{-}{\rm Mod}).$$
This completes the proof of $(1)$.

\vskip 5pt

$(2)$ \ By [Be, Theorem 4.10] we have that every module in $A\mbox{-}\mathcal{GP}$ is a filtered
colimit of modules in $A\mbox{-}{\mathcal G}proj$. This follows that $C_{gpac}(A\mbox{-}{\rm Mod})$ is closed under filtered colimits. By
minor modifications of the proof of [K2, Lemma 2$(4)$] we get that  $C_{gpac}(A\mbox{-}{\rm Mod})$ is closed under $\alpha$-pure subobjects for some
regular cardinal $\alpha$. Thus by [K2, Theorem 4] we get that every complex in $K(A\mbox{-}{\rm Mod})$ admits a right $K_{gpac}(A\mbox{-}{\rm Mod})$-approximation. Applying [N5, Proposition 1.4], it follows that the inclusion $i: K_{gpac}(A\mbox{-}$ ${\rm  Mod})\to K(A\mbox{-}{\rm Mod})$ admits a right adjoint $r$. Therefore we obtain from [Mi, Theorem 2.2] the following right recollement
$$
 K_{gpac}(A\mbox{-}{\rm Mod})\begin{smallmatrix}
 \underrightarrow{ \ \ \ i \ \ \ } \\
  \overleftarrow{ \ \ \ r \ \ \ }
\end{smallmatrix}
K(A\mbox{-}{\rm Mod})
\begin{smallmatrix}
  \underrightarrow{\ \ \  \ \ \ } \\
  \overleftarrow{ \ \ \  \ \ \ }
\end{smallmatrix}
D_{gp}(A\mbox{-}{\rm Mod})
$$
By [GZ1, Proposition I.1.3] we know that the right adjoint of the  quotient functor $K(A\mbox{-}$ ${\rm Mod})\to D_{gp}(A\mbox{-}{\rm Mod})$ is fully faithful. Therefore $D_{gp}(A\mbox{-}{\rm Mod})$ has small Hom-sets.

\vskip 5pt

Since $A$ is virtually Gorenstein of finite CM-type, it follows from the proof of $(1)$ that the inclusion
$K(A\mbox{-}\mathcal{GP})\to K(A\mbox{-}{\rm Mod})$ admits a right adjoint. Moreover, if $A$ is  Gorenstein, then by [G1, Theorem 2.7] we get that $K(A\mbox{-}\mathcal{GP})^{\perp}=K_{gpac}(A\mbox{-}{\rm Mod})$. This means that the inclusion $K_{gpac}(A\mbox{-}{\rm Mod})\to K(A\mbox{-}{\rm Mod})$ admits a left adjoint and the composition $K(A\mbox{-}\mathcal{GP})\to K(A\mbox{-}{\rm Mod})\to D_{gp}(A\mbox{-}{\rm Mod})$ is a triangle-equivalence. This completes the proof of $(2)$.
\hfill$\blacksquare$
\vskip 10pt

\begin{cor} \ Let $A$ be a Gorenstein  algebra of finite CM-type. Then the canonical functor $D_{gp}(A\mbox{-}{\rm Mod})\to D(A\mbox{-}{\rm Mod})$ admits
left and right adjoints that are fully faithful. The left adjoint preserves compactness and its restriction to compact objects identifies with the inclusion
$K^{b}(A\mbox{-}proj)\to K^{b}(A\mbox{-}\mathcal{G}proj)$.
\end{cor}

\vskip 5pt

\noindent{\bf Proof.} \ Since $A$ is Gorenstein of finite CM-type, it follows from Theorem 3.2 that $K(A\mbox{-}{\mathcal GP})$ is compactly generated, and $K(A\mbox{-}{\mathcal GP})^{c}\cong K^{b}(A\mbox{-}\mathcal{G}proj)$. By Corollary 3.3$(2)$ we have a triangle-equivalence $D_{gp}(A\mbox{-}{\rm Mod})\cong K(A\mbox{-}{\mathcal GP})$.  This implies  that $D_{gp}(A\mbox{-}{\rm Mod})^{c}\cong K^{b}(A\mbox{-}\mathcal{G}proj)$. Note that the  canonical functor $F: D_{gp}(A\mbox{-}{\rm Mod})\to D(A\mbox{-}{\rm Mod})$  preserves set-indexed coproducts and products. Thus by [N2, Theorem 4.1] we get that $F$ admits a left adjoint and a right adjoint, also the left adjoint preserves compactness. This follows from [N1, Theorem 2.1] that the restriction of  this left adjoint to compact objects identifies with the inclusion $K^{b}(A\mbox{-}proj)\to K^{b}(A\mbox{-}\mathcal{G}proj)$. Since $F$ is a quotient functor, it follows from [GZ1, Proposition I.1.3] that these adjoints are fully faithful.
\hfill$\blacksquare$
\vskip 10pt

Next we will compare the subcategory of compact objects in  the homotopy category of projective modules with the subcategory of compact objects in the homotopy category of Gorenstein projective modules.
We first recall the construction of compact generators of the homotopy category of projective $A$-modules for an artin algebra $A$ we will use.

\vskip 10pt

Let $A$ be an artin algebra  and $M$  a finitely generated $A$-module.
Then for $M^{*}$ there is a canonical quasi-isomorphism $P^{\bullet}\to M^{*}$ with $P^{\bullet}$ in  the homotopy category $K^{-,b}(A^{\rm op}\mbox{-}{\rm proj})$. We consider the collection of the form $P^{\bullet*}[n](n\in \mathbb{Z})$. Denote by $\mathcal{G}_{1}$ the set of one object from each such isomorphism class.

\vskip 10pt

\begin{lem} \ ([J, Theorem 2.4 and 3.2]) \ The homotopy category $K(A\mbox{-}\mathcal{P})$ is a compactly generated triangulated category with $\mathcal{G}_{1}$
as a set of compact generators. Moreover, there is a triangle-equivalence
$$K(A\mbox{-}\mathcal{P})^{c}\cong D^{b}(A^{\rm op})^{\rm op}$$
\end{lem}

\vskip 10pt

\begin{thm} \ Let $A$ be a virtually Gorenstein  algebra of finite CM-type. Then we have

\vskip 5pt

$(1)$ \ $K(A\mbox{-}\mathcal{GP})/K(A\mbox{-}\mathcal{P})$ is a compactly generated triangulated category.

\vskip 5pt

$(2)$ \ the localisation sequence of triangulated categories
$$K(A\mbox{-}\mathcal{P})\stackrel{i}\hookrightarrow K(A\mbox{-}\mathcal{GP})\stackrel{q}\to K(A\mbox{-}\mathcal{GP})/K(A\mbox{-}\mathcal{P})$$
yields, by restriction to compact objects, a sequence of functors
$$D^{b}(A^{\rm op})^{\rm op}\to D^{b}_{gp}(A^{\rm op})^{\rm op}\to (K(A\mbox{-}\mathcal{GP})/K(A\mbox{-}\mathcal{P}))^{c}.$$
Moreover the induced functor
$$F: D^{b}_{gp}(A^{\rm op})^{\rm op}/D^{b}(A^{\rm op})^{\rm op}\to (K(A\mbox{-}\mathcal{GP})/K(A\mbox{-}\mathcal{P}))^{c}$$
is fully faithful, and identifies $D^{b}_{gp}(A^{\rm op})^{\rm op}/D^{b}(A^{\rm op})^{\rm op}$ with a subcategory of
$(K(A\mbox{-}\mathcal{GP})$ $/K(A\mbox{-}\mathcal{P}))^{c}$ whose \'{e}paisse closure is all of $(K(A\mbox{-}\mathcal{GP})/K(A\mbox{-}\mathcal{P}))^{c}$.
\end{thm}

\vskip 5pt

\noindent{\bf Proof.} \ [G2, Theorem 2.6] implies $(1)$. Now we prove  $(2)$. By Lemma 3.5 we know that $K(A\mbox{-}\mathcal{P})$ is a compactly generated triangulated category with $\mathcal{G}_{1}$
as a set of compact generators, and there is an equivalence of triangulated categories $K(A\mbox{-}\mathcal{P})^{c}\cong D^{b}(A^{\rm op})^{\rm op}$.
By Theorem 3.2 we get that $K(A\mbox{-}\mathcal{GP})$ is a compactly generated triangulated category with $\mathcal{G}_{2}$ as a set of compact generators, and there is an equivalence of triangulated categories $K(A\mbox{-}\mathcal{GP})^{c}\cong D^{b}_{gp}(A^{\rm op})^{\rm op}$.

\vskip 5pt

Let $\{G^{\bullet}_{i}\}_{i\in I}$ be any family objects in  $K(A\mbox{-}\mathcal{GP})$. Then  ${\rm Hom}_{K(A\mbox{-}\mathcal{GP})}(iP, \coprod_{i\in I}G^{\bullet}_{i})={\rm Hom}_{K(A\mbox{-}\mathcal{GP})}(P, \coprod_{i\in I}G^{\bullet}_{i})\cong \coprod_{i\in I}{\rm Hom}_{K(A\mbox{-}\mathcal{GP})}(P, G^{\bullet}_{i})=\coprod_{i\in I}{\rm Hom}_{K(A\mbox{-}\mathcal{GP})}(iP, G^{\bullet}_{i})$ for each module $P\in A\mbox{-}{\rm proj}$.
So by [CFH, Proposition 2.6] the inclusion $i: K(A\mbox{-}\mathcal{P})\hookrightarrow K(A\mbox{-}\mathcal{GP})$ preserves compact objects.

\vskip 5pt

Applying [N1, Theorem 2.1] to the homotopy category $K(A\mbox{-}\mathcal{GP})$ and $K(A\mbox{-}\mathcal{P})$, we get that $i$ carries $K(A\mbox{-}\mathcal{P})^{c}$ to $K(A\mbox{-}\mathcal{GP})^{c}$, $q$ carries $K(A\mbox{-}\mathcal{GP})^{c}$ to $(K(A\mbox{-}\mathcal{GP})/K(A\mbox{-}\mathcal{P}))^{c}$, the natural functor $\widetilde{F}: K(A\mbox{-}\mathcal{GP})^{c}/K(A\mbox{-}\mathcal{P})^{c}
\to (K(A\mbox{-}\mathcal{GP})/K(A\mbox{-}\mathcal{P}))^{c}$ is fully faithful, and any object in $K(A\mbox{-}\mathcal{GP})^{c}/K(A\mbox{-}\mathcal{P})^{c}$
is a direct summand of some object in $(K(A\mbox{-}\mathcal{GP})/K(A\mbox{-}\mathcal{P}))^{c}$. This completes the proof of $(2)$.
\hfill$\blacksquare$
\vskip 10pt

Recall that an additive category $\mathcal{C}$ is  idempotent-complete if every idempotent morphism splits. Any additive category admits an idempotent
completion $l: \mathcal{C}\to \mathcal{C}^{\natural}$. Moreover, if $\mathcal{C}$ is triangulated, then $\mathcal{C}^{\natural}$ inherits a unique structure of triangulated category such that $l$ is a triangle functor(see [BS]).

\vskip 10pt

\begin{rem} \ Let $A$ be a virtually Gorenstein  algebra of finite CM-type.  Theorem 3.6 implies that the idempotent
completion of $D^{b}_{gp}(A^{\rm op})^{\rm op}/D^{b}(A^{\rm op})^{\rm op}$ and $(K(A\mbox{-}\mathcal{GP})/K(A\mbox{-}\mathcal{P}))^{c}$ are triangluar equivalent, also  $D^{b}(A)$  can be viewed as a triangulated subcategory of  $D^{b}_{gp}(A)$. But in general, this doesn't hold.
\end{rem}

\vskip 10pt

In general, a Gorenstein flat module is not necessarily Gorenstein projective for a ring.
Next we will establish a pair of  adjoint functors between the homotopy category of its Gorenstein flat modules and  the homotopy category of
a nice ring.

\vskip 10pt

\begin{thm} \ Let $R$ be a two-sided noetherian ring such that $R\mbox{-}\mathcal{GF}$ is closed under direct products. Then the inclusion $K(R\mbox{-}\mathcal{GF})\to K(R\mbox{-}{\rm Mod})$ has a right adjoint.
\end{thm}

\vskip 5pt

\noindent{\bf Proof.} \ By [EEI, Theorem 4.3] we get that every complex $X^{\bullet}$ in $K(R\mbox{-}{\rm Mod})$ admits  a $K(R\mbox{-}\mathcal{GF})$-precover.
Hence by [N5, Proposition 1.4] we get that the inclusion $K(R\mbox{-}\mathcal{GF})\to K(R\mbox{-}{\rm Mod})$ has a right adjoint.
\hfill$\blacksquare$
\vskip 10pt

Now we show an interesting phenomenon. Let $M$ be any $R$-module, and consider the complex below
$$\cdots\to 0\to M\to 0\to \cdots$$
The existence of a right adjoint to the inclusion gives us a morphism $Z^{\bullet}\to M$, in the category $K(R\mbox{-}{\rm mod})$,
\[\xymatrix{
&\cdots\ar [r] & Z^{-1} \ar [r]\ar[d]& Z^{0} \ar[d]^{\rho} \ar[r] &
Z^{1} \ar[d] \ar[r] &\cdots
 \\
& \cdots \ar [r] & 0 \ar [r] &
M \ar[r] & 0 \ar[r] & \cdots}\]
where the complex $Z^{\bullet}$ is a complex of Gorenstein flat $R$-modules. Furthermore, given any map $\varphi: F^{\bullet}\to M$,
with $F$ a Gorenstein flat $R$-module, we have a factorization of $\varphi:$
\[\xymatrix{
&\cdots\ar [r] & 0 \ar [r]\ar[d]& F \ar[d] \ar[r] &
0 \ar[d] \ar[r] &\cdots
 \\
&\cdots\ar [r] & Z^{-1} \ar [r]\ar[d]& Z^{0} \ar[d]^{\rho} \ar[r] &
Z^{1} \ar[d] \ar[r] &\cdots
 \\
& \cdots \ar [r] & 0 \ar [r] &
M \ar[r] & 0 \ar[r] & \cdots}\]
This shows the map $\rho: Z^{0}\to M$ is a Gorenstein flat precover for $M$.

\vskip 10pt

\begin{cor} \ Let $R$ be a two-sided noetherian ring such that $R\mbox{-}\mathcal{GF}$ is closed under direct products. Then the inclusion $i: K(R\mbox{-}\mathcal{F})\to K(R\mbox{-}\mathcal{GF})$ has a right adjoint.
\end{cor}

\vskip 5pt

\noindent{\bf Proof.} \ By [N5, Theorem 3.2] we get that the inclusion $i_{1}: K(R\mbox{-}\mathcal{F})\to K(R\mbox{-}{\rm Mod})$ has a right adjoint $j_{1}$. By Theorem 3.8 we get the inclusion $i_{2}: K(R\mbox{-}\mathcal{GF})\to K(R\mbox{-}{\rm Mod})$ has a right adjoint $j_{2}$.
Since we have a series of isomorphisms for any $F^{\bullet}\in K(R\mbox{-}\mathcal{F})$ and $Z^{\bullet}\in K(R\mbox{-}\mathcal{GF})$
\begin{align*}
{\rm Hom}_{K(R\mbox{-}\mathcal{GF})}(iF^{\bullet}, Z^{\bullet})&\cong {\rm Hom}_{K(R\mbox{-}{\rm Mod})}(i_{2}iF^{\bullet}, i_{2} Z^{\bullet})\\&=
{\rm Hom}_{K(R\mbox{-}{\rm Mod})}(i_{1}F^{\bullet}, i_{2} Z^{\bullet})\\&\cong {\rm Hom}_{K(R\mbox{-}\mathcal{F})}(F^{\bullet}, j_{1}i_{2} Z^{\bullet}),
\end{align*}
it follows that $i: K(R\mbox{-}\mathcal{F})\to K(R\mbox{-}\mathcal{GF})$ admits a right adjoint $j_{1}i_{2}$.
\hfill$\blacksquare$

\vskip 10pt

Enochs and Estrada in [EE] defined the category of  quasi-coherent $\Re$-modules, $\Re\mbox{-}{\rm mod}$, where $\Re$ is a representation by rings of a quiver $Q$. They aim to understand the category of quasi-coherent sheaves $\mathcal{O}coX$ on a scheme $X$ via the equivalence between it and $\Re\mbox{-}{\rm mod}$ for some quiver $Q$ and ring $\Re$. They also found that if $X$ is a locally Gorenstein  scheme, then $\mathcal{O}coX$ is a  Gorenstein category.
Now, an example of above theorem arise.
\vskip 10pt

\begin{exm} \ Let $A$ be a commutative noetherian ring and $(X, \mathcal{O}_{X})\subseteq \mathbb{P}^{n}(A)$ be a locally Gorenstein scheme. Denote by $\mathcal{Q}coX$ the category of quasi-coherent sheaves on $X$ and by $\mathcal{GF}X$ the subcategory of Gorenstein flat quasi-coherent
$\mathcal{O}_{X}$-modules.
Then the inclusion $K(\mathcal{GF}X)\to K(\mathcal{O} coX)$ has a right adjoint.
\end{exm}

\vskip 5pt

\noindent{\bf Proof.} \ Since $X$ is a locally Gorenstein scheme, following the notations in Section 3 in [EEG-R], $\mathcal{O} coX$ is equivalent to the category
$\Re\mbox{-}{\rm mod}$ as abelian categories, where $\Re$ is an associated ring such that $\Re(v)$ is a commutative Gorenstein ring for any vertex $v$.

\vskip 5pt

By [EX, Lemma 3.5] we know that $\Re(v)$-module $N$ is Gorenstein flat if and only if $N^{+}={\rm Hom}_{\mathbb{Z}}(N, \mathbb{Q}/\mathbb{Z})$ is a Gorenstein injective $\Re(v)$-module for all vertex $v$. So by [EEG-R, Corollary 3.13] we get that $M$  is a Gorenstein flat $\Re$-module if and only if $M(v)$ is a Gorenstein flat $\Re(v)$-module for all vertex $v$. Since $\Re(v)$ has a dualizing complex for all vertex $v$, it follows from [CFH, Theorem 5.7] that
$\Re(v)\mbox{-}\mathcal{GF}$ is closed under direct products. This deduces that $\Re\mbox{-}\mathcal{GF}$ is closed under direct products. By Theorem 3.8
we get  the inclusion $K(\Re\mbox{-}\mathcal{GF})\to K(\Re\mbox{-}{\rm Mod})$ admits a right adjoint. This deduces that $K(\mathcal{GF}X)\to K(\mathcal{O} coX)$ has a right adjoint.
\hfill$\blacksquare$

\vskip 20pt

\section{\bf Gorenstein representation Dimension}

\vskip 10pt

In this section we introduce the notion of Gorenstein representation dimension for an algebra of finite CM-type, and establish  relations among the dimension of its relative Auslander algebra,  Gorenstein representation dimension, the dimension of
the bounded Gorenstein derived category, and also the dimension of the bounded homotopy category of Gorenstein projective modules.

\vskip 10pt

Let $\mathcal{T}$ be a triangulated category, and $M\in \mathcal{T}$. We set
$$\langle M\rangle=\langle M\rangle_{1}={\rm add}\{M[i]|i\in\mathbb{Z}\}$$
$$\langle M\rangle_{n+1}={\rm add}\{X\mid \exists M^{'}\to X\to M^{''}\to M^{'}[1] \ with \ M^{'}\in \langle M\rangle, M^{''}\in \langle M\rangle_{n}\}$$
Recall from [Ro] that the dimension of a triangulated category $\mathcal{T}$ is the number
$${\rm dim}\mathcal{T}={\rm inf}\{n\in \mathbb{N}|\ there \ exists \ a \ M\in \mathcal{T} \ with \ \langle M\rangle_{n+1}=\mathcal{T}\}$$
For a subcategory $\mathcal{C}\subseteq \mathcal{T}$ the dimension is defined to be
${\rm dim}_{\mathcal{T}}\mathcal{C}={\rm inf}\{n\mid\exists M\in \mathcal{T}:\mathcal{C}\subseteq\langle M\rangle_{n+1}\}$.
We first have

\vskip 10pt

\begin{lem} \  Let $A$ be a  finite dimensional $k$-algebra of finite CM-type over a field $k$ and ${\mathcal G}p{\rm dim}X=n$ for some $X\in A\mbox{-}{\rm mod}$. Then $X\notin \langle G\rangle_{n}$.
\end{lem}

\vskip 5pt

\noindent{\bf Proof.}\quad The proper $A\mbox{-}{\mathcal G}proj$-resolution
$$\Omega^{n}_{G}X\rightarrowtail G_{n-1}\to G_{n-2}\to\cdots\to G_{1}\to G_{0}\twoheadrightarrow X$$
gives rise to a sequence of maps
$$X\to \Omega_{G}X[1]\to \cdots\to \Omega^{n-1}_{G}X[n-1]\to \Omega^{n}_{G}X[n]$$
in $D^{b}_{gp}(A)$. They are all $\langle G\rangle\mbox{-}$ghosts, and their composition in non-zero. Hence the claim follows from the ghost lemma
in the sense of Rouquier [Ro].
\hfill$\blacksquare$

\vskip 10pt

So we try to introduce the notion of Gorenstein representation dimension. Of particular interest to us
are ${\rm dim}{D^{b}_{gp}(A)}$ and ${\rm dim}_{D^{b}_{gp}(A)}(A\mbox{-}{\rm mod})$.

\vskip 10pt

\begin{defn} \ Let $A$ be an artin algebra of finite CM-type. The Gorenstein representation dimension of $A$ is defined as
$${\rm Grepdim}A={\rm min}\{M\mbox{-}{\rm resol.dim}(A\mbox{-}{\rm mod})|M\in A\mbox{-}{\rm mod} \ such \ that$$ $$G\oplus \nu(G)\in {\rm add}M\}+2.$$
An $A\mbox{-}$module $M$ realizing the minimum above is called the Gorenstein Auslander generator.
\end{defn}

\vskip 10pt

\begin{rem} \ Let $A$ be an artin algebra of finite CM-type. Then  ${\rm repdim}A\leq{\rm Grepdim}A$ by the definition of representation dimension (see [Au]). Also, if $A$ is CM-free (i.e. $A\mbox{-}\mathcal{GP}=A\mbox{-}\mathcal{P})$, then these two definitions are coincide.
\end{rem}

\vskip 10pt

\begin{lem} \ Let $A$ be an artin algebra of finite CM-type. Then
$${\rm Grepdim}A={\rm min}\{{\rm gl.dim}{\rm End}_{A}(M)\mid M\in A\mbox{-}{\rm mod} \ such \ that \ G\oplus \nu(G)\in {\rm add}M\}$$
\end{lem}

\vskip 5pt

\noindent{\bf Proof.}\quad For any $M\in A\mbox{-}{\rm mod}$ with $G\oplus \nu(G)\in {\rm add}M$,  we easily see that $M$ is a generator and cogenerator. Hence the claim follows from [Au]. \hfill$\blacksquare$

\vskip 10pt

Now we will establish relations among  above-mentioned dimensions for an algebra of finite CM-type.

\vskip 10pt

\begin{thm}\ Let $A$ be a  finite dimensional $k$-algebra of finite CM-type. Let $M\in A\mbox{-}{\rm mod}$ be a Gorenstein  Auslander generator, and $X\in A\mbox{-}{\rm mod}$. Then for any $n\in \mathbb{Z}$ we have
$$M\mbox{-}{\rm resol.dim}X\leq n\Longrightarrow X\in \langle M\rangle_{n+1}.$$
In particular,
$${\rm Grepdim}A\geq {\rm dim}_{D^{b}_{gp}(A)}(A\mbox{-}{\rm mod})+2.$$
\end{thm}

\vskip 5pt

\noindent{\bf Proof.}\quad This follows immediately from the fact that $A\mbox{-}\mathcal{G}proj\mbox{-}$acyclic short exact sequences in $A\mbox{-}{\rm mod}$ turn into distinguished triangles in $D^{b}_{gp}(A)$.
\hfill$\blacksquare$

\vskip 10pt

\begin{thm} \ Let $A$ be a  finite dimensional $k$-algebra of finite CM-type. Then
$${\rm Grepdim}A\geq {\rm dim}D^{b}_{gp}(A)\geq {\rm dim}\mathcal{G}p(A)-1$$
\end{thm}

\vskip 5pt

\noindent{\bf Proof.}\quad Let $M$ be a Gorenstein Auslander generator. By induction we get that every bounded complex of $A\mbox{-}{\rm mod}$ is $A\mbox{-}\mathcal{G}proj\mbox{-}$quasi-isomorphic to a bounded complex of ${\rm add}M$, i.e., the canonical functor $K^{b}({\rm add}M)\to D^{b}_{gp}(A)$ is essentially
surjective. Notice that we have canonical equivalences $K^{b}({\rm add}M)\cong K^{b}({\rm End}_{A}(M)\mbox{-}{\rm proj})\cong D^{b}({\rm End}_{A}(M))$ and that ${\rm dim}D^{b}({\rm End}_{A}(M))\leq {\rm gl.dim}{\rm End}_{A}(M)$. These deduce  that
${\rm Grepdim}A\geq {\rm dim}D^{b}_{gp}(A)$.

\vskip 5pt

Since $A$ is of finite CM-type, we have  that $D^{b}_{gp}(A)$ is triangular equivalent to $K^{-, gpb}(A\mbox{-}$ $\mathcal{G}proj)$,
and so $G$ is a generator of $D^{b}_{gp}(A)$. Note from [BS1, Theorem 2.8] that $D^{b}_{gp}(A)$ is idempotent split, and that ${\rm Hom}_{D^{b}_{gp}(A)}(G, -):  D^{b}_{gp}(A)\to {\rm Ab}$ is a cohomological functor such that ${\rm Hom}_{D^{b}_{gp}(A)}^{*}(G, G)=\mathcal{G}p(A)$ is an artin $k$-algebra. Hence  by [BIKO, Theorem 4.5] we get that ${\rm dim}D^{b}_{gp}(A)\geq {\rm dim}\mathcal{G}p(A)- 1$. This completes the proof.
\hfill$\blacksquare$

\vskip 10pt

\begin{thm} \ Let $A$ be a virtually Gorenstein  algebra of finite CM-type. Then we have ${\rm dim}K(A\mbox{-}\mathcal{GP})\geq {\rm dim}\mathcal{G}p(A)-1$.
\end{thm}

\vskip 5pt

\noindent{\bf Proof.} \  Since $A$ is virtually Gorenstein of finite CM-type, it follows from [Be, Theorem 4.10] that $G$ is a generator of $K(A\mbox{-}\mathcal{GP})$. Note that $K(A\mbox{-}\mathcal{GP})$ is idempotent split and ${\rm Hom}_{K(A\mbox{-}\mathcal{GP})}(G, -): K(A\mbox{-}\mathcal{GP})\to {\rm Ab}$ is a cohomological functor such that ${\rm Hom}_{K(A\mbox{-}\mathcal{GP})}^{*}(G, G)=\mathcal{G}p(A)$ is an artin $k$-algebra. Hence  by [BIKO, Theorem 4.5] we get that ${\rm dim}K(A\mbox{-}\mathcal{GP})\geq {\rm dim}\mathcal{G}p(A)-1$.
\hfill$\blacksquare$
\vskip 10pt

\vskip 20pt

{\bf Acknowledgements.} Some of this work was done when the author is visiting Professor Steffen K\"{o}nig in the University of Stuttgart, Germany.
The author  would like to thank  Steffen K\"{o}nig  for useful discussions and comments related to this work.

\end{document}